\documentclass{amsart}
\usepackage{amssymb}

\usepackage{amscd}

\newtheorem{theorem}{Theorem}[section]
\newtheorem{corollary}[theorem]{Corollary}

\newtheorem{lemma}[theorem]{Lemma}

\theoremstyle{definition}

\theoremstyle{remark}

\begin{document}
\title{A Petri theorem for rank-$2$ vector bundles with canonical determinant}
\author{Elisa Casini}
\author{Herbert Clemens}
\date{June, 2002}
\email{casini@science.unitn.it\\
clemens@math.utah.edu}
\maketitle

\begin{abstract}
This paper establishes the correctness of a conjecture of Bertram-Feinberg
and Mukai for a special class of globally generated rank-two bundles with
canonical determinant over a generic Riemann surface of genus at least four.
\end{abstract}

\section{Introduction}

\subsection{The problem}

Petri's general conjecture establishes the unobstructedness of linear series
on a general compact Riemann surface $C$ of genus $g>1$. (See, for example, 
\cite{C2}.) Some years ago, Bertram-Feinberg \cite{BF} and Mukai \cite{M}
independently formulated an analogous conjecture for stable rank-$2$ vector
bundles on $C$ with determinant $\omega _{C}$. The conjecture is that the
natural map 
\begin{equation}
Sym^{2}H^{0}\left( E\right) \rightarrow H^{0}\left( Sym^{2}E\right)
\label{weak}
\end{equation}
is injective for \textit{all} stable $E$. The stronger assertion that 
\begin{equation}
H^{0}\left( E\right) \otimes H^{0}\left( E\right) \rightarrow H^{0}\left(
E\otimes E\right)  \label{strong}
\end{equation}
is injective is known to be false for some stable bundles \cite{T}. In this
paper, over a generic Riemann surface $C$ of genus $g>3$, we construct
examples of globally generated, semi-stable bundles $E$ for which the
stronger assertion is true. We accomplish this by employing Hitchin's theory
of spectral curves \cite{H} to reduce the problem to an assertion about line
bundles $L$ of a spectral cover 
\begin{equation*}
\tilde{C}\rightarrow C
\end{equation*}
and then applying to $\left( \tilde{C},L\right) $ the techniques used in 
\cite{C2} to establish Petri's general conjecture. A critical ingredient
will be a condition introduced by Beauville \cite{B} which equates, under
sufficiently general conditions, the local analytic deformation space of $E$
with twists of $L$ by skew-symmetric line bundles on $\tilde{C}$. The
necessity of this condition is shown by an example of V. Mercat, which we
explore in an Appendix to this paper.

\subsection{Hitchin's spectral curves}

Let $E$ be a stable rank-$2$ vector bundle over $C$ with 
\begin{equation*}
\det E=\omega ,
\end{equation*}
the canonical line bundle on $C$. We have a natural isomorphism 
\begin{eqnarray}
E &\rightarrow &Hom\left( E,\omega \right)  \label{Niso} \\
e &\mapsto &e\wedge  \notag
\end{eqnarray}
Following the theory of Hitchin \cite{H}, let 
\begin{equation*}
\tilde{\varphi}\in \mathrm{Hom}\left( E^{\vee },E\right) =\mathrm{Hom}\left(
E^{\vee },E^{\vee }\otimes \omega \right) .
\end{equation*}
Replacing $\varphi $ by 
\begin{equation*}
\tilde{\varphi}-\frac{tr\tilde{\varphi}}{2}
\end{equation*}
if necessary, we can (and will) assume that 
\begin{equation}
\tilde{\varphi}\in \Gamma \left( Hom^{0}\left( E^{\vee },E^{\vee }\otimes
\omega \right) \right)  \label{phitilde}
\end{equation}
where $Hom^{0}$ means homomorphisms of trace zero. We assume that $\tilde{%
\varphi}$ is not the zero homomorphism. Now 
\begin{equation*}
y^{2}+\varphi
\end{equation*}
is the characteristic polynomial of $\tilde{\varphi}$ where 
\begin{equation*}
\varphi =\det \tilde{\varphi}\in H^{0}\left( \omega ^{2}\right) .
\end{equation*}
and $y$ takes values in sections of $\omega $. Considering $\varphi \left(
C\right) $ as a curve in the geometric line bundle $\omega ^{2}$, let $%
\tilde{C}$ denote the inverse image of $C$ under the squaring map 
\begin{equation*}
\omega \rightarrow \omega ^{2}.
\end{equation*}
We then have a double covering 
\begin{equation*}
\pi :\tilde{C}\rightarrow C
\end{equation*}
branched at 
\begin{equation*}
\Delta =zero\left( \varphi \right) ,
\end{equation*}
and the arithmetic genus $\tilde{g}$ of $\tilde{C}$ is computed from the
identity 
\begin{eqnarray*}
2\tilde{g}-2 &=&4\left( 2g-2\right) \\
\tilde{g} &=&2\left( 2g-2\right) +1.
\end{eqnarray*}

\subsection{Globally generated $E$}

As in \S 4 of \cite{vGI}, suppose that $E$ is globally generated. We can
then some construct morphisms $\left( \ref{phitilde}\right) $ as follows.
Let 
\begin{equation}
W\subseteq H^{0}\left( E\right)  \label{W}
\end{equation}
be any subspace such that the evaluation map 
\begin{equation*}
W\otimes \mathcal{O}_{C}\rightarrow E
\end{equation*}
is surjective. Then we have the standard exact sequence 
\begin{equation}
0\rightarrow F\rightarrow W\otimes \mathcal{O}_{C}\rightarrow E\rightarrow 0.
\label{Wquot}
\end{equation}
So for each quadric 
\begin{equation*}
Q\in Sym^{2}W
\end{equation*}
we obtain a diagram 
\begin{equation*}
\begin{array}{ccccccccc}
0 & \rightarrow & E^{\vee } & \overset{\varepsilon ^{\vee }}{\longrightarrow 
} & W^{\vee }\otimes \mathcal{O}_{C} & \rightarrow & F^{\vee } & \rightarrow
& 0 \\ 
&  &  &  & \downarrow ^{T} &  &  &  &  \\ 
0 & \rightarrow & F & \rightarrow & W\otimes \mathcal{O}_{C} & \overset{
\varepsilon }{\longrightarrow } & E & \rightarrow & 0
\end{array}
\end{equation*}
where 
\begin{equation*}
T\left( w^{\vee }\right) =\left\langle \left. Q\right| w^{\vee
}\right\rangle \in W^{\vee \vee }=W.
\end{equation*}
So we have an induced morphism 
\begin{equation*}
\tilde{Q}=\varepsilon \circ T\circ \varepsilon ^{\vee }\in \Gamma
Hom^{0}\left( E^{\vee },E\right) .
\end{equation*}
Rewriting 
\begin{eqnarray*}
Hom^{0}\left( E^{\vee },E\right) &=&Hom^{0}\left( E^{\vee },E^{\vee }\otimes
\omega \right) =Sym^{2}E, \\
\tilde{Q} &\mapsto &\left( \left( e_{1}^{\vee },e_{2}^{\vee }\right) \mapsto
\left\langle \left. \tilde{Q}\left( e_{1}^{\vee }\right) \right| e_{2}^{\vee
}\right\rangle \right)
\end{eqnarray*}
the map 
\begin{eqnarray*}
Sym^{2}W &\rightarrow &H^{0}\left( Sym^{2}E\right) \\
Q &\mapsto &\tilde{Q}
\end{eqnarray*}
is just the standard map induced by multiplication of sections.

On the other hand, for any section $y$ of $\omega $, we have 
\begin{equation}
y\cdot :E^{\vee }\rightarrow E^{\vee }\otimes \omega =E  \label{ydot}
\end{equation}
with trace $2y.$ For fixed $c\in C$, let 
\begin{equation*}
e_{1},e_{2}
\end{equation*}
be a basis for $E_{c}$ and let 
\begin{equation*}
e_{1}\wedge e_{2}=y_{c}\in \omega _{c}.
\end{equation*}
Then the mapping $\left( \ref{ydot}\right) $ is given by 
\begin{eqnarray*}
e_{1}^{\vee } &\mapsto &\frac{-y\left( c\right) }{y_{c}}e_{2} \\
e_{2}^{\vee } &\mapsto &\frac{y\left( c\right) }{y_{c}}e_{1}.
\end{eqnarray*}
and $\tilde{Q}_{c}$ is given by 
\begin{eqnarray}
e_{1}^{\vee } &\mapsto &Q_{11}\cdot e_{1}+Q_{12}\cdot e_{2}  \label{quadric}
\\
e_{2}^{\vee } &\mapsto &Q_{21}\cdot e_{1}+Q_{22}\cdot e_{2}.  \notag
\end{eqnarray}
where 
\begin{equation*}
Q_{ij}=Q\left( e_{i}^{\vee },e_{j}^{\vee }\right) .
\end{equation*}
The mapping 
\begin{equation*}
\left( y_{c}\tilde{Q}_{c}-y\left( c\right) \cdot \right) :E_{c}^{\vee
}\rightarrow E_{c}
\end{equation*}
drops rank if and only 
\begin{equation*}
y\left( c\right) =\pm y_{c}\sqrt{\det \left( 
\begin{array}{cc}
Q_{11} & Q_{12} \\ 
Q_{21} & Q_{22}
\end{array}
\right) }
\end{equation*}
with eigenvectors given by the zeros of the quadratic form 
\begin{equation*}
\left. Q\right| _{E_{c}^{\vee }}\text{.}
\end{equation*}

\subsection{$E$ as a pushforward of a line bundle}

For smooth $\tilde{C}$, the eigenvector subspace 
\begin{equation*}
\left\{ \left( c,v\right) \in E^{\vee }:\tilde{\varphi}\left( v\right)
=\left( \pm \sqrt{\varphi }\right) v\right\}
\end{equation*}
has the structure of a line bundle $L^{\vee }$ on $\tilde{C}$. The
restriction map 
\begin{equation*}
Hom\left( E^{\vee },\mathcal{O}_{C}\right) \rightarrow Hom\left( L^{\vee },%
\mathcal{O}_{C}\right) =:L
\end{equation*}
induces an isomorphism 
\begin{equation*}
E=\pi _{*}L.
\end{equation*}
Since 
\begin{equation*}
\det \pi _{*}L=\omega ,
\end{equation*}
we have 
\begin{equation*}
\deg L=2\left( 2g-2\right) =\tilde{g}-1.
\end{equation*}
For any small deformation $L^{\prime }$ of $\left\{ L\right\} $ in $Pic^{%
\tilde{g}-1}\left( \tilde{C}\right) $, $E^{\prime }=\pi _{*}L^{\prime }$ is
stable, but 
\begin{equation*}
\det E^{\prime }
\end{equation*}
may not equal $\omega $, a difficulty we can always remedy by twisting $%
L^{\prime }$ by $\pi ^{*}\left( M\right) $ for some $\left\{ M\right\} \in
Pic^{0}\left( C\right) $.

Let 
\begin{equation*}
\begin{array}{r}
\iota :\tilde{C}\rightarrow \tilde{C} \\ 
y\mapsto -y
\end{array}
\end{equation*}
denote the involution on $\tilde{C}$. By considering $\pm 1$-eigenspaces, $%
\iota $ induces an isomorphism 
\begin{equation}
\pi _{*}\mathcal{O}_{\tilde{C}}=\mathcal{O}_{C}\oplus \mathcal{O}_{C}\left( 
\sqrt{-\Delta }\right) =\mathcal{O}_{C}\oplus \mathcal{O}_{C}\left( \omega
^{-1}\right) =\mathcal{O}_{C}\oplus T_{C}.  \label{k1}
\end{equation}
A first-order deformation of $L$ is given by an element of 
\begin{equation*}
H^{1}\left( \mathcal{O}_{\tilde{C}}\right) =H^{1}\left( \mathcal{O}%
_{C}\right) \oplus H^{1}\left( \mathcal{O}_{C}\left( \sqrt{-\Delta }\right)
\right) .
\end{equation*}
By the projection formula, first-order deformations in $H^{1}\left( \mathcal{%
O}_{C}\right) $ change $\det \pi _{*}L$. Thus first-order deformations with
fixed determinant are given by elements of 
\begin{equation*}
H^{1}\left( \mathcal{O}_{C}\left( \sqrt{-\Delta }\right) \right) \subseteq
H^{1}\left( \mathcal{O}_{\tilde{C}}\right) .
\end{equation*}
The natural map 
\begin{equation*}
\pi _{*}End\left( L\right) \rightarrow End\left( \pi _{*}L\right)
\end{equation*}
is an element of 
\begin{equation*}
\mathcal{O}_{C}\oplus Hom\left( \mathcal{O}_{C}\left( \sqrt{-\Delta }\right)
,End^{0}E\right)
\end{equation*}
where $End^{0}$ denotes those endomorphism of trace zero. The summand of
this map in 
\begin{equation*}
Hom\left( \mathcal{O}_{C}\left( \sqrt{-\Delta }\right) ,End^{0}E\right)
=Hom^{0}\left( E^{\vee },E\right)
\end{equation*}
is given by 
\begin{equation*}
\tilde{\varphi}\in \Gamma \left( Hom^{0}\left( E^{\vee },E^{\vee }\otimes
\omega \right) \right)
\end{equation*}
considered as the inclusion 
\begin{equation}
\tilde{\varphi}:\mathcal{O}_{C}\left( \sqrt{-\Delta }\right) \rightarrow
End^{0}E  \label{incl0}
\end{equation}
The induced mapping 
\begin{equation}
H^{1}\left( \mathcal{O}_{C}\left( \sqrt{-\Delta }\right) \right) \rightarrow
H^{1}\left( End^{0}E\right)  \label{iso1}
\end{equation}
is an isomorphism for general choice of $\left( C,\tilde{\varphi}\right) $.
(See \cite{B}, \S 1.5.)

\subsection{$\tilde{\varphi}$ arising from a quadric}

Returning to the exact sequence $\left( \ref{Wquot}\right) $, the dual
sequence 
\begin{equation*}
0\rightarrow E^{\vee }\rightarrow W^{\vee }\otimes \mathcal{O}
_{C}\rightarrow F^{\vee }\rightarrow 0
\end{equation*}
lets us view $\left( \ref{Wquot}\right) $ as a morphism 
\begin{eqnarray*}
l &:&C\rightarrow G:=Gr\left( 2,W^{\vee }\right) . \\
c &\mapsto &\Bbb{P}\left( E_{c}^{\vee }\right)
\end{eqnarray*}
So we see that, if we choose a quadric $Q$ on $\Bbb{P}\left( W^{\vee
}\right) $, then for 
\begin{equation}
\tilde{\varphi}=\tilde{Q}  \label{goodcase}
\end{equation}
we have that $\tilde{C}$ is the inverse image of $C$ under the double cover
of $\Bbb{P}\left( W\right) $ branched at $c\in C$ exactly when the line $%
l\left( c\right) $ is tangent to the quadric $Q$. So$,$ if we choose $Q$
generally, $\tilde{C}$ has $2g-2$ distinct branchpoints in $C$ and so must
be smooth. (The authors wish to thank Christian Pauly for pointing this fact
out to us and providing the above proof.) Thus $\tilde{C}$ is the double
cover of $C$ induced from the standard double cover of $G$ induced by the
quadric $Q$, that is, the double cover branched along the divisor consisting
of those lines tangent to $Q$.

Continuing in the situation $\left( \ref{goodcase}\right) $, a section 
\begin{equation*}
w\in W\subseteq H^{0}\left( E\right)
\end{equation*}
gives a section of $L$ which vanishes exactly when one of the two points of 
\begin{equation*}
\Bbb{P}\left( E_{c}^{\vee }\right) \cap Q
\end{equation*}
lies in the hyperplane $w\in W=W^{\vee \vee }$. So 
\begin{equation}
\pi _{*}L=\mathcal{O}_{G}\left( 2\right) .  \label{Dtilde}
\end{equation}

\subsection{Quotients of $C\times \Bbb{C}^{3}$}

If $E$ is globally generated, we can choose $W$ in $\left( \ref{W}\right) $
such that 
\begin{equation}
\dim W=3.  \label{good}
\end{equation}
In this case, we obtain a morphism 
\begin{eqnarray}
C &\rightarrow &\Bbb{P}\left( W\right) =\Bbb{P}^{2}.  \label{Cmap} \\
c &\mapsto &\Bbb{P}\left( F_{c}\right)  \notag
\end{eqnarray}
Since 
\begin{equation*}
c_{1}\left( F\right) =\omega ^{-1}
\end{equation*}
we see that $\left( \ref{Cmap}\right) $ is a projection of the canonical
mapping$.$ Conversely suppose we take any projection of the canonical curve 
\begin{equation*}
p:C\rightarrow \Bbb{P}^{2},
\end{equation*}
whose center does not meet the canonical curve. Then the resulting exact
sequence 
\begin{equation*}
0\rightarrow p^{*}\mathcal{O}_{\Bbb{P}^{2}}\left( -1\right) \overset{p}{%
\longrightarrow }\Bbb{C}^{3}\otimes \mathcal{O}_{C}\rightarrow E\rightarrow 0
\end{equation*}
gives rise to a globally generated rank-$2$ bundle $E$ with canonical
determinant. If a section of $E$ coming from $w\in \Bbb{C}^{3}$ has $m$
zeros, then 
\begin{equation*}
p\left( C\right)
\end{equation*}
must have an $m$-fold point at $w$ and conversely$.$ Assuming $C$ is not
hyperelliptic, by a result of Lazarsfeld (see Theorem 1.1 of \cite{G}) there
are three sections $\alpha ,\beta ,\gamma \in H^{0}\left( \omega \right) $
such that the map 
\begin{equation*}
\alpha \cdot +\beta \cdot +\gamma \cdot :H^{0}\left( \omega \right)
\rightarrow H^{0}\left( \omega ^{2}\right)
\end{equation*}
is surjective. Dually, for general $C$, the generic map 
\begin{equation*}
H^{1}\left( T_{C}\right) \rightarrow H^{1}\left( \mathcal{O}_{C}\right)
^{\oplus 3}
\end{equation*}
induced by the choice of $3$ sections of $\omega $ is injective, so that 
\begin{equation}
\Bbb{C}^{3}\rightarrow H^{0}\left( E\right)  \label{surjtoE}
\end{equation}
is surjective for generic $C$ and $p$. Now any non-semi-stable bundle $E$
has a sub-bundle of degree at least $g$ and therefore a section with $g$ or
more zeros. Therefore $p\left( C\right) $ has a $g$-tuple point. Thus we
conclude that, for general choice of $C$ and $p$, the bundle $E$ is
semi-stable if $g>2$ since $p\left( C\right) $ has only nodes. Furthermore,
for general choice of $C$ and $p$ and generic 
\begin{equation*}
\tilde{\varphi}\in \Gamma Hom^{0}\left( E^{\vee },E\right) .
\end{equation*}
the associated spectral curve $\tilde{C}$ is smooth.

We also claim that the map $\left( \ref{iso1}\right) $ is an isomorphism in
this case. To see this last assertion, it will suffice to prove that $\left( 
\ref{iso1}\right) $ is an isomorphism for some \textit{special} choice of $%
C,p,\tilde{\varphi}$. For example let $C^{\prime }$ be a general curve which
has a vanishing theta-null, that is, admits a line bundle $J$ with $%
J^{2}=\omega $ and 
\begin{equation*}
h^{0}\left( J\right) =2.
\end{equation*}
(This is a codimension-one condition on $C$.) Let 
\begin{equation*}
E^{\prime }=J\oplus J.
\end{equation*}
Since $J$ is globally generated, so is $E^{\prime }$. If $g>3$, then $%
C^{\prime }$ is not hyperelliptic, so again by the result of Lazarsfeld
there are three sections $\alpha ,\beta ,\gamma \in H^{0}\left( \omega
\right) $ such that the map 
\begin{equation*}
\alpha \cdot +\beta \cdot +\gamma \cdot :H^{0}\left( \omega \right)
\rightarrow H^{0}\left( \omega ^{2}\right)
\end{equation*}
is surjective. Define 
\begin{equation*}
\tilde{\varphi}^{\prime }:J^{\vee }\oplus J^{\vee }\rightarrow J\oplus J
\end{equation*}
by the matrix 
\begin{equation*}
\left( 
\begin{array}{cc}
\alpha \cdot & \beta \cdot \\ 
\gamma \cdot & -\alpha \cdot
\end{array}
\right) .
\end{equation*}
The mapping 
\begin{eqnarray*}
\Gamma Hom^{0}\left( \left( E^{\prime }\right) ^{\vee },E^{\prime }\right)
&\rightarrow &H^{0}\left( \omega ^{2}\right) \\
\tilde{\psi} &\mapsto &Tr\left( \tilde{\psi}\circ \tilde{\varphi}\right)
\end{eqnarray*}
is therefore surjective. So by Serre duality, the map 
\begin{equation*}
\tilde{\varphi}^{\prime }:H^{1}\left( T_{C}\right) \rightarrow H^{1}\left(
End^{0}E^{\prime }\right)
\end{equation*}
is injective. (Compare \cite{B}, \S 1.3.) Now for generic 
\begin{equation*}
W\subseteq H^{0}\left( E^{\prime }\right) ,\ \dim W=3
\end{equation*}
$E^{\prime }$ is a quotient as in $\left( \ref{Wquot}\right) $ and so given
by a sub-bundle $F^{\prime }$ corresponding to a projection of the canonical
curve $C^{\prime }$ into $\Bbb{P}^{2}$. But the triple $\left( C^{\prime
},F^{\prime },\tilde{\varphi}^{\prime }\right) $ can be deformed to the
generic such triple $\left( C,F,\tilde{\varphi}\right) $, that is, to the
generic projection of a general curve $C$ of genus $g$ to $\Bbb{P}^{2}$. The
resulting quotient $E$ is globally generated by construction and $\left( \ref
{iso1}\right) $ for $\tilde{\varphi}$ is injective by semi-continuity of
rank. Thus $\tilde{\varphi}$ is an isomorphism since the general quotient $E$
has $h^{0}\left( End^{0}E\right) =0$

Finally in the situation $\left( \ref{good}\right) $ and $\left( \ref
{goodcase}\right) $, given $c\in C$, we can choose three conics $%
Q_{1},Q_{2},Q_{3}$ such that 
\begin{equation*}
\left. Q_{i}\right| _{E_{c}^{\vee }}
\end{equation*}
form a basis for $Sym^{2}E_{c}$. Thus for $\tilde{\varphi}_{i}=\tilde{Q}_{i}$
in $\left( \ref{incl0}\right) $, the map 
\begin{equation}
\tilde{\varphi}_{1}+\tilde{\varphi}_{2}+\tilde{\varphi}_{3}:\mathcal{O}
_{C}\left( \sqrt{-\Delta }\right) ^{\oplus 3}\rightarrow End^{0}E
\label{gensurj}
\end{equation}
is generically surjective. Thus the same is true for \textit{generic} choice
of $\tilde{\varphi}_{i}\in \Gamma Hom^{0}\left( E^{\vee },E\right) $.

\subsection{The theorem}

In what follows we shall treat only those rank-$2$ bundles $E$ such that:

\begin{description}
\item[Condition 1]  $E$ is simple, globally generated, $\det E=\omega $.

\item[Condition 2]  The map $\left( \ref{iso1}\right) $ is an isomorphism
for generic $\tilde{\varphi}\in \Gamma Hom^{0}\left( E^{\vee },E\right) $.
\end{description}

Notice that, since $\mathcal{O}_{C}\left( \sqrt{-\Delta }\right) =T_{C}$,
Condition 2) implies that 
\begin{equation*}
h^{0}\left( End^{0}E\right) =3g-3+0+3\left( 1-g\right) =0,
\end{equation*}
that is, $E$ is \textit{simple}. Since $H^{0}\left( E\right) =H^{0}\left(
L\right) $, global generation in Condition 1) requires that 
\begin{equation*}
h^{0}\left( L\right) \geq 3.
\end{equation*}
So we are restricting attention to a locally closed subvariety of relatively
high codimension in $\frak{M}\left( C\right) $.

Let 
\begin{equation*}
W_{2,\omega }^{r}\subseteq \frak{M}\left( C\right)
\end{equation*}
denote the scheme defined by those $E$ such that $h^{0}\left( E\right) \geq
r+1$. Bertram-Feinberg \cite{BF} and Mukai \cite{M} show that $W_{2,\omega
}^{r}$ has the natural structure of a determinantal scheme whose cotangent
space at $\left( C,E\right) $ is the cokernel of the map 
\begin{equation}
Sym^{2}H^{0}\left( E\right) \rightarrow H^{0}\left( Sym^{2}E\right)
\label{symmap}
\end{equation}
and that, if $\left( \ref{symmap}\right) $ is injective, then $W_{2,\omega
}^{r}$ is smooth and reduced and so of dimension 
\begin{equation*}
\rho \left( 2,r,\omega \right) =3g-3-\frac{\left( r+1\right) \left(
r+2\right) }{2}
\end{equation*}
at $\left( C,E\right) $. The purpose of this paper is to prove:

\begin{theorem}
\label{MT}If $C$ is of general moduli and $E$ satisfies Condition 1 and
Condition 2 above, the map 
\begin{equation*}
H^{0}\left( E\right) \otimes H^{0}\left( E\right) \rightarrow H^{0}\left(
E\otimes E\right) 
\end{equation*}
is injective, so that $\left( \ref{symmap}\right) $ is also. If $g\geq 4$,
there exist bundles $E$ satisfying Conditions $1$ and $2$ so that, by
semi-continuity, these conditions are satisfied for a generic globally
generated $E$ with canonical determinant.
\end{theorem}

For $C$ is a curve of genus $g$ of general moduli, this result has been
conjectured for \textit{all} stable rank-$2$ $E$ with canonical determinant
by Bertram-Feinberg and by Mukai (\textit{loc. cit.}). A corollary in the
special case in which $E$ satisfies Condition 1 and Condition 2 above is
that the map 
\begin{equation*}
\bigwedge\nolimits^{2}H^{0}\left( E\right) \rightarrow H^{0}\left(
\bigwedge\nolimits^{2}E\right) =H^{0}\left( \omega \right)
\end{equation*}
is also injective. So, in particular 
\begin{equation*}
3\leq \frac{h^{0}\left( E\right) \left( h^{0}\left( E\right) -1\right) }{2}%
\leq g.
\end{equation*}

\subsection{Analogy with proof of Petri conjecture}

Our proof will develop analogously to the proof of the classical Petri
conjecture for line bundles in \cite{C2}. In particular we employ Kuranishi
theory as in \cite{C1} and \cite{C2} to relate $n$-th order deformation
theory of $\left( C,\iota ,L\right) $ to the first cohomology of sheaves of
(holomorphic) differential operators 
\begin{equation*}
\frak{D}_{n}\left( L\right)
\end{equation*}
of order $\leq n$ on sections of $L$.

\section{The space $\Bbb{P}\left( E^{\vee }\right) $}

Let 
\begin{equation*}
X=\Bbb{P}\left( E^{\vee }\right) .
\end{equation*}
The natural inclusion 
\begin{equation*}
L^{\vee }\subseteq E^{\vee }
\end{equation*}
induces an imbedding 
\begin{equation}
\tilde{C}\rightarrow X  \label{incl1}
\end{equation}
which assigns to an eigenvalue $\left( c,\sqrt{\varphi \left( c\right) }%
\right) $ in $\tilde{C}$ the corresponding eigenvector in $E_{c}^{\vee }$.
Abusing notation we will let 
\begin{equation*}
\pi :X\rightarrow C
\end{equation*}
also denote the projection map extending $\pi :\tilde{C}\rightarrow C$.

Thus we have the exact sequences 
\begin{equation}
0\rightarrow T_{\pi }\rightarrow T_{X}\rightarrow \pi ^{*}T_{C}\rightarrow 0
\label{tanXseq}
\end{equation}
and 
\begin{equation*}
0\rightarrow \pi _{*}T_{\pi }\rightarrow \pi _{*}T_{X}\rightarrow
T_{C}\rightarrow 0.
\end{equation*}
Notice that the Euler sequence 
\begin{equation*}
0\rightarrow \mathcal{O}_{X}\rightarrow \pi ^{*}E^{\vee }\otimes \mathcal{O}
_{X}\left( 1\right) \rightarrow T_{\pi }\rightarrow 0
\end{equation*}
gives 
\begin{equation*}
T_{\pi }=\pi ^{*}T_{C}\otimes \mathcal{O}_{X}\left( 2\right)
\end{equation*}
and 
\begin{eqnarray*}
\pi _{*}T_{\pi } &=&End^{0}\left( E^{\vee }\right) =End^{0}E \\
&=&T_{C}\otimes Sym^{2}E.
\end{eqnarray*}
Since 
\begin{eqnarray*}
Ext^{1}\left( \pi ^{*}T_{C},T_{\pi }\right) &=&H^{1}\left( \pi ^{*}\omega
\otimes T_{\pi }\right) \\
&=&H^{1}\left( \omega \otimes \pi _{*}T_{\pi }\right) \\
&=&H^{1}\left( \omega \otimes End^{0}E\right) \\
&=&H^{0}\left( End^{0}E\right) ^{\vee } \\
&=&0
\end{eqnarray*}
the sequence $\left( \ref{tanXseq}\right) $ is split and hence so is its
push-forward. Thus non-canonically 
\begin{eqnarray}
T_{X} &=&T_{\pi }\oplus \pi ^{*}T_{C}  \label{split} \\
\pi _{*}T_{X} &=&End^{0}E\oplus T_{C}.
\end{eqnarray}

The Euler sequence also gives 
\begin{equation*}
0\rightarrow \mathcal{O}_{X}\rightarrow \pi ^{*}End\left( E^{\vee }\right)
\rightarrow \pi ^{*}\pi _{*}T_{\pi }\rightarrow 0.
\end{equation*}
The kernel of the mapping 
\begin{eqnarray*}
\pi ^{*}End\left( E^{\vee }\right) &\rightarrow &\pi ^{*}E^{\vee }\otimes 
\mathcal{O}_{X}\left( 1\right) \\
\left( c,\left[ v\right] ,\varepsilon \right) &\mapsto &\left( c,\left[
v\right] ,\left. \varepsilon \right| _{\left[ v\right] }\right)
\end{eqnarray*}
at $\left( c,\left[ v\right] \right) $ is those endomorphisms $\varepsilon $
which vanish on $v$ and so the kernel of 
\begin{equation}
\pi ^{*}End\left( E^{\vee }\right) \rightarrow T_{\pi }  \label{factor}
\end{equation}
is those $\varepsilon $ such that $v$ is an eigenvector of $\varepsilon $.
Thus the inclusion 
\begin{equation*}
\pi ^{*}T_{C}\rightarrow \pi ^{*}End\left( E^{\vee }\right)
\end{equation*}
in $\left( \ref{incl0}\right) $ given by 
\begin{equation*}
\tilde{\varphi}\in Hom\left( E^{\vee },E^{\vee }\otimes \omega \right)
\end{equation*}
composes with $\left( \ref{factor}\right) $ to give a map 
\begin{equation}
\pi ^{*}T_{C}\rightarrow T_{\pi }  \label{key}
\end{equation}
of line bundles on $X$ which vanishes exactly along $\tilde{C}$. Thus 
\begin{equation*}
T_{\pi }=\mathcal{O}_{X}\left( \tilde{C}\right) \otimes \pi ^{*}T_{C}.
\end{equation*}
By the Euler sequence 
\begin{equation*}
\pi ^{*}\det \left( E^{\vee }\right) \otimes \mathcal{O}_{X}\left( 2\right)
=T_{\pi }
\end{equation*}
so that 
\begin{equation}
\mathcal{O}_{X}\left( \tilde{C}\right) =\mathcal{O}_{X}\left( 2\right) .
\label{iso2}
\end{equation}

\section{The higher $\nu $-maps}

The isomorphism $\left( \ref{iso2}\right) $ and the resulting exact sequence 
\begin{equation*}
0\rightarrow \mathcal{O}_{X}\left( -1\right) \rightarrow \mathcal{O}%
_{X}\left( 1\right) \rightarrow L\rightarrow 0
\end{equation*}
gives, by applying $R\pi _{*}$, isomorphisms 
\begin{equation*}
H^{i}\left( E\right) =H^{i}\left( \mathcal{O}_{X}\left( 1\right) \right)
\rightarrow H^{i}\left( L\right)
\end{equation*}
for all $i\geq 0$. The tangent space to the deformation space of $\left( X,%
\mathcal{O}_{X}\left( 1\right) \right) $ is 
\begin{equation*}
H^{1}\left( \frak{D}_{1}\left( \mathcal{O}_{X}\left( 1\right) \right)
\right) =H^{1}\left( \pi _{*}\frak{D}_{1}\left( \mathcal{O}_{X}\left(
1\right) \right) \right)
\end{equation*}
and the element $\alpha \in H^{0}\left( L\right) =H^{0}\left( E\right) $
deforms to first-order under $\zeta \in H^{1}\left( \frak{D}_{1}\left( 
\mathcal{O}_{X}\left( 1\right) \right) \right) $ if and only if 
\begin{equation*}
\zeta \cdot \alpha \in H^{1}\left( \mathcal{O}_{X}\left( 1\right) \right)
=H^{1}\left( E\right)
\end{equation*}
vanishes. (See \cite{AC} or \cite{C2}.)

Next let 
\begin{equation*}
\frak{D}_{n}\subseteq \frak{D}_{n}\left( E\right)
\end{equation*}
denote the image of $\pi _{*}\frak{D}^{n}\left( \mathcal{O}_{X}\left(
1\right) \right) $ under the natural map 
\begin{equation*}
\pi _{*}\frak{D}_{n}\left( \mathcal{O}_{X}\left( 1\right) \right)
\rightarrow \frak{D}_{n}\left( \pi _{*}\mathcal{O}_{X}\left( 1\right)
\right) .
\end{equation*}
Using $\left( \ref{split}\right) $ we have the (symbol) exact sequences 
\begin{equation*}
0\rightarrow \mathcal{O}_{C}\rightarrow \frak{D}_{1}\overset{\sigma }{%
\longrightarrow }T_{C}\oplus End^{0}E\rightarrow 0
\end{equation*}
and 
\begin{equation*}
0\rightarrow \frak{D}_{n-1}\rightarrow \frak{D}_{n}\overset{\sigma }{%
\longrightarrow }T_{C}^{n}\oplus \left( T_{C}^{n-1}\otimes End^{0}E\right)
\rightarrow 0.
\end{equation*}
So we conclude from the symbol sequence 
\begin{equation*}
0\rightarrow \frak{D}_{n-1}\left( E\right) \rightarrow \frak{D}_{n}\left(
E\right) \overset{\sigma }{\longrightarrow }T_{C}^{n}\otimes EndE\rightarrow
0
\end{equation*}
that 
\begin{equation}
\frak{D}_{n-1}\left( E\right) \subseteq \frak{D}_{n}.  \label{indincl}
\end{equation}

From $\left( \ref{incl0}\right) $, we have the inclusion 
\begin{equation*}
T_{C}\oplus T_{C}\rightarrow T_{C}\oplus End^{0}E
\end{equation*}
induced by $\tilde{\varphi}$. We define 
\begin{equation*}
\frak{\tilde{D}}_{1}=\sigma ^{-1}\left( T_{C}\oplus T_{C}\right)
\end{equation*}
resulting in the (symbol) exact sequence 
\begin{equation*}
0\rightarrow \mathcal{O}_{C}\rightarrow \frak{\tilde{D}}_{1}\overset{\sigma 
}{\longrightarrow }T_{C}\oplus T_{C}\rightarrow 0.
\end{equation*}
Since 
\begin{equation*}
\frak{\tilde{D}}_{1}\subseteq \pi _{*}\frak{D}_{1}\left( \mathcal{O}
_{X}\left( 1\right) \right) \rightarrow \frak{D}_{1}\left( \pi _{*}\mathcal{O%
}_{X}\left( 1\right) \right)
\end{equation*}
is injective, we can (and will) consider $\frak{\tilde{D}}_{1}$ as a
subsheaf of $\frak{D}_{1}$. In fact 
\begin{equation*}
\frak{\tilde{D}}_{1}\cap \frak{D}_{0}\left( E\right) =\mathcal{O}_{C}\oplus
T_{C}\subseteq \frak{D}_{0}\left( E\right) =EndE
\end{equation*}
where the containment is given by sending $\mathcal{O}_{C}$ to scalar
endomorphisms and using the map 
\begin{equation*}
\tilde{\varphi}:T_{C}\rightarrow End^{0}E.
\end{equation*}
Also, $\left( \ref{iso1}\right) $ implies that 
\begin{equation*}
H^{1}\left( \frak{\tilde{D}}_{1}\right)
\end{equation*}
is the tangent space to the deformation space of pairs $\left( C,E\right) $.
The subspace 
\begin{equation*}
S_{\omega }\subseteq H^{1}\left( \frak{\tilde{D}}_{1}\right)
\end{equation*}
of those deformations such that $\det E$ remains the canonical bundle maps
isomorphically onto $H^{1}\left( T_{C}\oplus T_{C}\right) $.

Let 
\begin{equation*}
S=\left\{ \zeta \in S_{\omega }:\left( \zeta \cdot :H^{0}\left( E\right)
\rightarrow H^{1}\left( E\right) \right) =0\right\} .
\end{equation*}
The assumption that $C$ is generic implies that the map 
\begin{equation}
S\rightarrow H^{1}\left( \frac{\frak{\tilde{D}}_{1}}{\frak{\tilde{D}}%
_{1}\cap \frak{D}_{0}\left( E\right) }\right) =H^{1}\left( T_{C}\right)
\label{mu1}
\end{equation}
is surjective.

We let 
\begin{equation*}
\frak{\tilde{D}}_{n}\subseteq \frak{D}_{n}
\end{equation*}
be the image of the natural map 
\begin{equation*}
\frak{\tilde{D}}_{1}^{\otimes n}\rightarrow \frak{D}_{n}
\end{equation*}
induced by composition of operators. We have the morphism of exact sequences 
\begin{equation}
\begin{array}{ccccccccc}
0 & \rightarrow & \frak{\tilde{D}}_{n-1} & \rightarrow & \frak{\tilde{D}}_{n}
& \overset{\sigma }{\longrightarrow } & T_{C}^{n}\otimes Sym^{n}\left( 
\mathcal{O}_{C}^{\oplus 2}\right) & \rightarrow & 0 \\ 
&  & \downarrow &  & \downarrow &  & \downarrow &  &  \\ 
0 & \rightarrow & \frak{D}_{n-1} & \rightarrow & \frak{D}_{n} & \overset{
\sigma }{\longrightarrow } & T_{C}^{n}\oplus \left( T_{C}^{n-1}\otimes
End^{0}E\right) & \rightarrow & 0
\end{array}
\label{symbseq}
\end{equation}

As with the Petri problem for line bundles (see \cite{C2}), we want to
consider the maps 
\begin{equation*}
\tilde{\nu}^{n}:H^{1}\left( \frak{\tilde{D}}_{n}\right) \rightarrow \mathrm{%
\ Hom}\left( H^{0}\left( E\right) ,H^{1}\left( E\right) \right)
\end{equation*}
and the induced maps 
\begin{equation*}
\nu ^{n}:H^{1}\left( T_{C}^{n}\oplus \left( T_{C}^{n-1}\otimes Sym^{n}\left( 
\mathcal{O}_{C}^{\oplus 2}\right) \right) \right) \rightarrow \frac{\mathrm{%
\ Hom}\left( H^{0}\left( E\right) ,H^{1}\left( E\right) \right) }{image%
\tilde{ \nu}^{n-1}}.
\end{equation*}
The surjectivity of $\left( \ref{iso1}\right) $ shows that 
\begin{equation}
\tilde{\nu}^{1}\left( H^{1}\left( \frak{\tilde{D}}_{1}\right) \right)
=image\left( H^{1}\left( EndE\right) \rightarrow \mathrm{Hom}\left(
H^{0}\left( E\right) ,H^{1}\left( E\right) \right) \right) .  \label{step1}
\end{equation}

Analogously to the solution Petri problem for line bundles, we wish to show
that, for general $C$ and $E$ as above, 
\begin{equation}
\nu ^{n}=0  \label{n>1}
\end{equation}
for $n>1$. The next two sections are devoted to a proof of this fact. In the
next section, we construct elements of $H^{1}\left( \frak{\tilde{D}}
_{n}\right) $ whose symbols span $H^{1}\left( T_{C}^{n}\oplus \left(
T_{C}^{n-1}\otimes Sym^{n}\left( \mathcal{O}_{C}^{\oplus 2}\right) \right)
\right) $ for $n>1$. In the following section we show that these elements
lie in the kernel of $\tilde{\nu}^{n}$.

\section{Generators of the domains of the higher $\nu $-maps}

To construct the requisite elements of $H^{1}\left( \frak{\tilde{D}}%
_{n}\right) $, we need to choose divisors along which to make (Schiffer)
deformations of our general curve $C.$ First recall that 
\begin{eqnarray*}
H^{1}\left( \pi _{*}T_{X}\right) &=&H^{1}\left( T_{X}\right) \\
&=&H^{1}\left( T_{\pi }\oplus \pi ^{*}T_{C}\right) \\
&=&H^{1}\left( End^{0}E\right) \oplus H^{1}\left( T_{C}\right) \\
&=&H^{1}\left( T_{C}\right) \oplus H^{1}\left( T_{C}\right)
\end{eqnarray*}
is the tangent space to the deformations of $X$. The first factor
corresponds to deformations of the bundle (with canonical determinant) and
the second summand corresponds to first-order deformations of the base curve 
$C$. Suppose now that $A$ is a sufficiently ample simple divisor on $C$ and
let $U$ be a small anaytic neighborhood of (the support of) $A$. The
(non-canonical) splitting 
\begin{equation*}
T_{X}=T_{\pi }\oplus \pi ^{*}T_{C},
\end{equation*}
whose summands we will call the vertical and horizontal tangent spaces
respectively) allows us to fix an isomorphism 
\begin{equation*}
\pi ^{-1}\left( U\right) =U\times \Bbb{P}^{1}.
\end{equation*}
such that 
\begin{equation*}
\left. \pi ^{*}T_{C}\right| _{\pi ^{-1}\left( U\right) }=T_{U}\boxtimes 
\mathcal{O}_{\Bbb{P}^{1}}.
\end{equation*}
The inclusion 
\begin{equation*}
\pi ^{*}T_{C}\rightarrow T_{\pi }
\end{equation*}
which vanishes along $\tilde{C}$ also realizes $\pi ^{*}T_{C}$ as a sheaf of
(vertical) vector fields on $X$. Let 
\begin{equation*}
\alpha _{1},\beta _{1}
\end{equation*}
be $C^{\infty }$-sections of $T_{C}$ with support inside $U$ which are
meromorphic some neighborhood of $A$ with poles only above the points of $A$
. By the above remarks, we have distinguished liftings of $\alpha _{1},\beta
_{1}$ to vector fields on $X$ supported in $\pi ^{-1}\left( U\right) $ which
are respectively horizontal and vertical. We denote the lifted vector fields
again as $\alpha _{1},\beta _{1}$. Then, as in \cite{C1} and \cite{C2}, 
\begin{equation}
\left[ \overline{\partial },e^{-tL_{\alpha _{1}+\beta _{1}}}\right]
\label{formdef}
\end{equation}
determines a deformation of $X$. There is a lifting of $\alpha _{1}+\beta
_{1}$ to a vector field $\tilde{\alpha}_{1}+\tilde{\beta}_{1}$ on the total
space of 
\begin{equation*}
\mathcal{O}_{X}\left( -1\right)
\end{equation*}
with the property that 
\begin{equation*}
\left[ \chi ,\tilde{\alpha}_{1}\right] =0=\left[ \chi ,\tilde{\beta}%
_{1}\right]
\end{equation*}
where $\chi $ is the Euler vector field on $\mathcal{O}_{X}\left( -1\right) $
. Any two liftings differ by an element 
\begin{equation*}
\left( a\circ \pi \right) \chi
\end{equation*}
where $a$ is a function supported on $U$ and meromorphic near $A$ with poles
only above the points of $A$. The tangent space to $PicX$ is $H^{1}\left( 
\mathcal{O}_{C}\right) $ and therefore as in Lemma 3.2 of \cite{C2} we can
arrange that a deformation 
\begin{equation}
\left[ \overline{\partial },e^{-\left( tL_{\tilde{\alpha}_{1}+\tilde{\beta}%
_{1}+a_{1}\chi }+t^{2}L_{a_{2}\chi }+t^{3}L_{a_{3}\chi }+\ldots \right)
}\right]  \label{bunddef}
\end{equation}
of the bundle $\mathcal{O}_{X}\left( -1\right) $ has the property that 
\begin{equation*}
\det \pi _{*}\mathcal{O}_{X}\left( 1\right)
\end{equation*}
remains canonical.

We can modify $\left( \ref{bunddef}\right) $ over a given deformation of $C$
by modifying $\beta _{1}$ as follows. 
\begin{equation}
H^{0}\left( \pi ^{*}\frac{T_{C}\left( A\right) }{T_{C}}\right) \rightarrow
H^{1}\left( \pi ^{*}T_{C}\right) =H^{1}\left( End^{0}E\right)  \label{surj}
\end{equation}
is surjective. Recalling that $H^{1}\left( End^{0}E\right) $ is the tangent
space to $\frak{M}\left( C\right) $ at $\left\{ E\right\} $, we can
therefore choose a $C^{\infty }$-section $\varepsilon _{1}$ of $\pi
^{*}T_{C}\subseteq T_{\pi }$ which is supported over $U$ and meromorphic at $%
\pi ^{-1}\left( A\right) $ so that, applying $\overline{\partial }$ to the
element, we can realize any given element of $H^{1}\left( End^{0}E\right) $.
Thus 
\begin{equation}
\left[ \overline{\partial },e^{tL_{\alpha _{1}+\beta _{1}+\varepsilon
_{1}}}\right]  \label{modformdef}
\end{equation}
gives any first-order deformation of $X$ over the given deformation of $C$.
Again using Lemma 3.2 of \cite{C2} above we can arrange that a deformation 
\begin{equation}
\left[ \overline{\partial },e^{-\left( tL_{\tilde{\alpha}_{1}+\tilde{\beta}%
_{1}+\tilde{\varepsilon}_{1}+a_{1}\chi }+t^{2}L_{a_{2}\chi
}+t^{3}L_{a_{3}\chi }+\ldots \right) }\right]  \label{newmod}
\end{equation}
of the bundle $\mathcal{O}_{X}\left( -1\right) $ has the property that 
\begin{equation*}
\det \pi _{*}\mathcal{O}_{X}\left( 1\right)
\end{equation*}
remains canonical.

However the sections of 
\begin{equation*}
E=\pi _{*}\mathcal{O}_{X}\left( 1\right)
\end{equation*}
may not all deform under the deformation of $\left( X,\mathcal{O}_{X}\left(
1\right) \right) $ determined by $\left( \ref{newmod}\right) $. We will call
a deformation of $X$ \textit{admissible} if all sections of $\pi _{*}%
\mathcal{O}_{X}\left( 1\right) $ do indeed deform. To modify $\left( \ref
{modformdef}\right) $ to produce an admissible deformation, we proceed as
follows. Choose a sufficiently ample divisor $B$ supported in an open set $%
V\subseteq C$ disjoint from $U.$ Using the surjectivity of 
\begin{equation*}
H^{0}\left( \pi ^{*}\frac{T_{C}\left( B\right) }{T_{C}}\right) \rightarrow
H^{1}\left( \pi ^{*}T_{C}\right) =H^{1}\left( End^{0}E\right)
\end{equation*}
we can choose a $C^{\infty }$-section $\delta _{1}$ of $\pi ^{*}T_{C}$ which
is supported on $\pi ^{-1}\left( V\right) $ and meromorphic at $\pi
^{-1}\left( V\right) $ and, applying $\overline{\partial }$ to the element,
we can realize any given element of $H^{1}\left( End^{0}E\right) $. And we
can pick $\delta _{1}$ so that 
\begin{equation*}
\left[ \overline{\partial },e^{tL_{\alpha _{1}+\beta _{1}+\varepsilon
_{1}+\delta _{1}}}\right]
\end{equation*}
gives any first-order deformation of $X$ over the given deformation of $C$,
in particular, the admissible one which must exist by the assumption that $C$
is generic. Similarly we can choose $\delta _{2}$ so that 
\begin{equation*}
\left[ \overline{\partial },e^{tL_{\alpha _{1}+\beta _{1}+\varepsilon
_{1}+\delta _{1}}+t^{2}L_{\delta _{2}}}\right]
\end{equation*}
gives the admissible second-order deformation over the given deformation of $%
C$, etc. Repeating this argument we achieve an admissible formal deformation 
\begin{equation}
\left[ \overline{\partial },e^{tL_{\alpha _{1}+\beta _{1}+\varepsilon
_{1}+\delta _{1}}+t^{2}L_{\delta _{2}}+t^{3}L_{\delta _{3}}+\ldots }\right]
\label{adformdef}
\end{equation}
of $X$ over the deformation $\left( \ref{formdef}\right) $ of $C$. Again
using Lemma 3.2 of \cite{C2} above we can arrange that a deformation 
\begin{equation}
\left[ \overline{\partial },e^{-\left( tL_{\tilde{\alpha}_{1}+\tilde{\beta}%
_{1}+\tilde{\varepsilon}_{1}+\tilde{\delta}_{1}+a_{1}\chi
}+t^{2}L_{a_{2}\chi +\tilde{\delta}_{2}}+t^{3}L_{a_{3}\chi +\tilde{\delta}%
_{3}}+\ldots \right) }\right]  \label{liftdef}
\end{equation}
of the bundle $\mathcal{O}_{X}\left( -1\right) $ has the property that 
\begin{equation*}
\det \pi _{*}\mathcal{O}_{X}\left( 1\right)
\end{equation*}
remains canonical.

Now we are ready to look at the symbols of the operators we have just
constructed. In 
\begin{equation*}
H^{1}\left( T_{C}^{n}\right) \otimes H^{0}\left( Sym^{n}\left( \mathcal{O}
_{C}^{\oplus 2}\right) \right) =H^{1}\left( T_{C}^{n}\otimes Sym^{n}\left( 
\mathcal{O}_{C}^{\oplus 2}\right) \right)
\end{equation*}
we consider the subspace 
\begin{equation*}
H^{1}\left( T_{C}^{n}\right) \otimes H^{0}\left( Sym^{n}\left( \mathcal{O}
_{C}\oplus 0\right) \right)
\end{equation*}
corresponding to the $n$-th power of \textit{vertical} tangent vectors, that
is, those lying in $T_{\pi }.$

\begin{lemma}
\label{almost}i) Elements $\left( \ref{liftdef}\right) $ are of the form 
\begin{equation*}
\sum\nolimits_{n>0}D_{n}t^{j}
\end{equation*}
with 
\begin{equation*}
D_{n}\in H^{1}\left( \frak{\tilde{D}}_{n}\right) .
\end{equation*}

ii) For $n>1$, the maps $\nu _{+}^{n}$ factor through 
\begin{equation*}
H^{1}\left( \frac{T_{C}^{n}\otimes Sym^{n}\left( \mathcal{O}_{C}^{\oplus
2}\right) }{T_{C}^{n}\otimes Sym^{n}\left( \mathcal{O}_{C}\oplus 0\right) }%
\right) 
\end{equation*}

iii) For $n>1$, the collection of such elements $D_{n}$ obtained by varying $%
\tilde{\alpha}_{1},\tilde{\varepsilon}_{1}$ and the divisor $A$ have symbols 
\begin{equation*}
\sigma \left( D_{n}\right) 
\end{equation*}
which span $H^{1}\left( \frac{T_{C}^{n}\otimes Sym^{n}\left( \mathcal{O}%
_{C}^{\oplus 2}\right) }{T_{C}^{n}\otimes Sym^{n}\left( \mathcal{O}%
_{C}\oplus 0\right) }\right) $.
\end{lemma}

\begin{proof}
i) Consider the (direct limit) exact sequence 
\begin{equation*}
0\rightarrow \frak{\tilde{D}}_{n}\rightarrow \mathcal{O}_{C}\left( \infty
\cdot \left( A+B\right) \right) \otimes \frak{\tilde{D}}_{n}\rightarrow 
\frac{\mathcal{O}_{C}\left( \infty \cdot \left( A+B\right) \right) \otimes 
\frak{\tilde{D}}_{n}}{\frak{\tilde{D}}_{n}}\rightarrow 0.
\end{equation*}
The coefficient of $t^{n}$ in $e^{-\left( tL_{\tilde{\alpha}_{1}+\tilde{\beta%
}_{1}+\tilde{\varepsilon}_{1}+\tilde{\delta}_{1}+a_{1}\chi
}+t^{2}L_{a_{2}\chi +\tilde{\delta}_{2}}+t^{3}L_{a_{3}\chi +\tilde{\delta}%
_{3}}+\ldots \right) }$ lies in $H^{0}\left( \frac{\mathcal{O}_{C}\left(
\infty \cdot \left( A+B\right) \right) \otimes \frak{\tilde{D}}_{n}}{\frak{%
\tilde{D}}_{n}}\right) $.

ii) Suppose first that we restrict the elements $\left( \ref{liftdef}\right) 
$ to those for which $\alpha _{1}=0$. Then the symbols of the coefficients
of $t^{n}$ lie in the ``subspace of vertical symbols'' 
\begin{equation*}
H^{1}\left( T_{C}^{n}\right) \otimes H^{0}\left( Sym^{n}\left( \mathcal{O}
_{C}\oplus 0\right) \right) \subseteq H^{1}\left( T_{C}^{n}\right) \otimes
H^{0}\left( Sym^{n}\left( \mathcal{O}_{C}^{\oplus 2}\right) \right) .
\end{equation*}
Letting the subscript ``$n$'' denote the coefficient of $t^{n}$, we have by
Lemma 4.1 of \cite{C2} that we can vary $\varepsilon _{1}$ and $A$ to obtain
elements 
\begin{equation*}
D_{n}=\left[ \overline{\partial },e^{-\left( tL_{0+\tilde{\beta}_{1}+\tilde{%
\varepsilon}_{1}+\tilde{\delta}_{1}+a_{1}\chi }+t^{2}L_{a_{2}\chi +\tilde{%
\delta}_{2}}+t^{3}L_{a_{3}\chi +\tilde{\delta}_{3}}+\ldots \right) }\right]
_{n}
\end{equation*}
whose symbols generate 
\begin{equation*}
H^{1}\left( T_{C}^{n}\right) \otimes H^{0}\left( Sym^{n}\left( \mathcal{O}
_{C}\oplus 0\right) \right) \subseteq H^{1}\left( T_{C}^{n}\right) \otimes
H^{0}\left( Sym^{n}\left( \mathcal{O}_{C}^{\oplus 2}\right) \right) .
\end{equation*}
Furthermore these elements $D_{n}$ have the property that 
\begin{equation*}
\tilde{\nu}^{n}\left( D_{n}\right) \in image\left( H^{1}\left(
End^{0}E\right) \rightarrow Hom\left( H^{0}\left( E\right) ,H^{1}\left(
E\right) \right) \right) \subseteq image\left( \nu ^{1}\right) .
\end{equation*}
Thus for $n>1$ we can consider $\nu ^{n}$ as a map 
\begin{equation*}
\nu ^{n}:H^{1}\left( \frac{T_{C}^{n}\otimes Sym^{n}\left( \mathcal{O}%
_{C}^{\oplus 2}\right) }{T_{C}^{n}\otimes Sym^{n}\left( \mathcal{O}%
_{C}\oplus 0\right) }\right) \rightarrow \frac{Hom\left( H^{0}\left(
E\right) ,H^{1}\left( E\right) \right) }{image\left( \nu ^{n-1}\right) }.
\end{equation*}
iii) Notice that our restriction on the support of the $\delta _{j}$ insures
that, for $n>1$, the coefficient of $t^{n}$ in $\left( \ref{liftdef}\right) $
has the \textit{same} symbol in 
\begin{equation*}
H^{1}\left( \frac{T_{C}^{n}\otimes Sym^{n}\left( \mathcal{O}_{C}^{\oplus
2}\right) }{T_{C}^{n}\otimes Sym^{n}\left( \mathcal{O}_{C}\oplus 0\right) }%
\right)
\end{equation*}
as 
\begin{equation}
\left\{ \left[ \overline{\partial },e^{tL_{\alpha _{1}+\beta
_{1}+\varepsilon }}\right] \right\}  \label{twsym}
\end{equation}
does. So it suffices to show that the symbols of the coefficients of $t^{n}$
in elements of the form $\left( \ref{twsym}\right) $ generate 
\begin{equation*}
H^{1}\left( T_{C}^{n}\right) \otimes H^{0}\left( Sym^{n}\left( \mathcal{O}%
_{C}^{\oplus 2}\right) \right) =H^{1}\left( T_{C}^{n}\otimes Sym^{n}\left( 
\mathcal{O}_{C}^{\oplus 2}\right) \right) .
\end{equation*}
But $Sym^{n}H^{0}\left( \mathcal{O}_{C}^{\oplus 2}\right) $ generates $%
H^{0}\left( Sym^{n}\left( \mathcal{O}_{C}^{\oplus 2}\right) \right) $, so it
suffices to show that elements of the form 
\begin{equation*}
\left[ \overline{\partial },\tilde{\alpha}_{1}^{n}\right]
\end{equation*}
generate $H^{1}\left( T_{C}^{n}\right) $. But this is again the content of
Lemma 4.1 of \cite{C2}.
\end{proof}

\section{Vanishing higher $\nu $-maps}

\begin{lemma}
\label{crit}For generic $C$ and for $E$ satisfying Condition 1 and Condition
2, the maps 
\begin{equation*}
\nu ^{n}:H^{1}\left( T_{C}^{n}\oplus \left( T_{C}^{n-1}\otimes Sym^{n}\left( 
\mathcal{O}_{C}^{\oplus 2}\right) \right) \right) \rightarrow \frac{\mathrm{%
Hom}\left( H^{0}\left( E\right) ,H^{1}\left( E\right) \right) }{image\left( 
\tilde{\nu}^{n-1}\right) }
\end{equation*}
vanish for $n>1$.
\end{lemma}

\begin{proof}
By Lemma \ref{almost} therefore, it suffices to show that, for each
expression $\left( \ref{liftdef}\right) $, there are operators 
\begin{equation*}
\tilde{D}_{n}\in H^{1}\left( \frak{\tilde{D}}_{n}\right)
\end{equation*}
such that

1) $\tilde{D}_{n}$ has the same symbol in 
\begin{equation*}
H^{1}\left( T_{C}^{n}\otimes Sym^{n}\left( \mathcal{O}_{C}^{\oplus 2}\right)
\right)
\end{equation*}
as the coefficient of $t^{n}$ in 
\begin{equation}
\left[ \overline{\partial },e^{-\left( tL_{\tilde{\alpha}_{1}+\tilde{\beta}
_{1}+\tilde{\varepsilon}_{1}+\tilde{\delta}_{1}+a_{1}\chi
}+t^{2}L_{a_{2}\chi +\tilde{\delta}_{2}}+t^{3}L_{a_{3}\chi +\tilde{\delta}
_{3}}+\ldots \right) }\right] ,  \label{Kdata}
\end{equation}

2) 
\begin{equation*}
\tilde{\nu}^{n}\left( \tilde{D}_{n}\right) =0.
\end{equation*}
The proof of this fact is identical to the proof given in \S 3.4-3.5 of \cite
{C2}. The notational dictionary is 
\begin{eqnarray*}
reference\ paper &\leftrightarrow &this\ paper \\
X_{0} &\leftrightarrow &X \\
L_{0} &\leftrightarrow &\mathcal{O}_{X}\left( 1\right)  \\
S^{n}\left( T_{X_{0}}\right)  &\leftrightarrow &T_{C}^{n}\otimes
Sym^{n}\left( \mathcal{O}_{C}^{\oplus 2}\right) .
\end{eqnarray*}
In the proof we consider the line bundle 
\begin{equation*}
\tilde{E}\left( 1\right) =\mathcal{O}_{\Bbb{P}\left( H^{0}\left( E\right)
\right) }\left( 1\right) \boxtimes \mathcal{O}_{X}\left( 1\right) 
\end{equation*}
on the product 
\begin{equation*}
\Bbb{P}\left( H^{0}\left( \mathcal{O}_{X}\left( 1\right) \right) \right) 
\times X.
\end{equation*}
Just as in \S 3.4-3.5 of \cite{C2} we use Lemma 2.9 and Lemma 3.1iii) of 
\textit{loc. cit.} to twist the Kuranishi data $\left( \ref{Kdata}\right) $
above to obtain 
\begin{equation}
\left[ \overline{\partial },e^{-\left( tL_{\tilde{\alpha}_{1}+\tilde{\beta}%
_{1}+\tilde{\varepsilon}_{1}+\tilde{\delta}_{1}+a_{1}\chi
}+t^{2}L_{a_{2}\chi +\tilde{\delta}_{2}}+t^{3}L_{a_{3}\chi +\tilde{\delta}%
_{3}}+\ldots \right) }e^{-L_{\tilde{\beta}_{T}\left( t\right) }}\right] 
\label{K'data}
\end{equation}
for which the condition that the tautological section extends under
deformation implies that the map 
\begin{equation*}
H^{0}\left( \tilde{E}\left( 1\right) \right) \rightarrow
\sum\nolimits_{n>0}H^{1}\left( \tilde{E}\left( 1\right) \right) t^{n}
\end{equation*}
given by $\left( \ref{K'data}\right) $ is $0$.
\end{proof}

Next notice that the generic surjectivity of $\left( \ref{gensurj}\right) $
implies that the maps 
\begin{equation*}
\tilde{\varphi}_{1}+\tilde{\varphi}_{2}+\tilde{\varphi}_{3}:\left( \mathcal{O%
}_{C}\left( \sqrt{-\Delta }\right) ^{\oplus 3}\right) \otimes
T_{C}^{n-1}\rightarrow \left( End^{0}E\right) \otimes T_{C}^{n-1}
\end{equation*}
are all generically surjective. Thus the maps 
\begin{equation}
H^{1}\left( \left( \mathcal{O}_{C}\left( \sqrt{-\Delta }\right) \right)
\otimes T_{C}^{n-1}\right) ^{\oplus 3}=H^{1}\left( \left( \mathcal{O}%
_{C}\left( \sqrt{-\Delta }\right) ^{\oplus 3}\right) \otimes
T_{C}^{n-1}\right) \rightarrow H^{1}\left( \left( End^{0}E\right) \otimes
T_{C}^{n-1}\right)  \label{triplesurj}
\end{equation}
are all surjective. For $i=1,2,3$, let 
\begin{equation*}
\frak{\tilde{D}}_{n}^{\left( i\right) }
\end{equation*}
be the sheaf $\frak{\tilde{D}}_{n}$ on $C$ derived from the spectral curve
associated to $\tilde{\varphi}_{i}$. Now define a map 
\begin{equation*}
F:\frak{\tilde{D}}_{n}^{\left( 1\right) }\oplus \frak{\tilde{D}}_{n}^{\left(
2\right) }\oplus \frak{\tilde{D}}_{n}^{\left( 3\right) }\rightarrow \frak{D}
_{n}
\end{equation*}
by adding the three maps 
\begin{equation*}
\frak{\tilde{D}}_{n}^{\left( i\right) }\rightarrow \frak{D}_{n}
\end{equation*}
associated to the spectral curves for $\tilde{\varphi}_{i}$ for $i=1,2,3$.
We have a commutative diagram 
\begin{equation*}
\begin{array}{ccccccccc}
0 & \rightarrow & \sum\nolimits_{i}\frak{\tilde{D}}_{n-1}^{\left( i\right) }
& \rightarrow & \sum\nolimits_{i}\frak{\tilde{D}}_{n}^{\left( i\right) } & 
\rightarrow & \left( T_{C}^{n}\otimes Sym^{n}\left( \mathcal{O}_{C}^{\oplus
2}\right) \right) ^{\oplus 3} & \rightarrow & 0 \\ 
&  & \downarrow ^{F} &  & \downarrow ^{F} &  & \downarrow ^{F} &  &  \\ 
0 & \rightarrow & \frak{D}_{n-1} & \rightarrow & \frak{D}_{n} & \rightarrow
& T_{C}^{n}\oplus \left( \left( End^{0}E\right) \otimes T_{C}^{n-1}\right) & 
\rightarrow & 0
\end{array}
\end{equation*}
and the generic surjectivity of $\left( \ref{triplesurj}\right) $ and
induction shows that 
\begin{equation}
F_{*}:H^{1}\left( \sum\nolimits_{i}\frak{\tilde{D}}_{n}^{\left( i\right)
}\right) \rightarrow H^{1}\left( \frak{D}_{n}\right)  \label{lastsurj}
\end{equation}
is surjective for each $n$.

\begin{corollary}
Letting 
\begin{equation*}
\tilde{\mu}^{n}:H^{1}\left( \frak{D}_{n}\left( E\right) \right) \rightarrow 
\mathrm{Hom}\left( H^{0}\left( E\right) ,H^{1}\left( E\right) \right) 
\end{equation*}
denote the natural maps induced by the action of $\frak{D}_{n}\left(
E\right) $ on sections of $E$, the induced maps 
\begin{equation*}
\mu ^{n}:H^{1}\left( T_{C}^{n}\otimes EndE\right) \rightarrow \frac{\mathrm{%
\ Hom}\left( H^{0}\left( E\right) ,H^{1}\left( E\right) \right) }{%
image\left( \tilde{\mu}^{n-1}\right) }
\end{equation*}
are all zero for $n\geq 1$.
\end{corollary}

\begin{proof}
Referring to $\left( \ref{indincl}\right) $, consider the chain 
\begin{equation*}
\frak{D}_{0}=\mathcal{O}_{C}\subseteq \frak{D}_{0}\left( E\right)
=EndE\subseteq \frak{D}_{1}\subseteq \frak{D}_{1}\left( E\right) \subseteq 
\frak{D}_{2}\subseteq \frak{D}_{2}\left( E\right) \subseteq \ldots
\end{equation*}
of subsheaves and the image of $H^{1}$ of each in $\mathrm{Hom}\left(
H^{0}\left( E\right) ,H^{1}\left( E\right) \right) $. By $\left( \ref{step1}
\right) $, Lemma \ref{crit}, and the surjectivity of $\left( \ref{lastsurj}
\right) $, the induced maps on $H^{1}$ of the graded quotients of this chain
into respective quotients of $\mathrm{Hom}\left( H^{0}\left( E\right)
,H^{1}\left( E\right) \right) $ are all zero from 
\begin{equation*}
H^{1}\left( \frac{\frak{D}_{1}}{EndE}\right) \rightarrow \frac{\mathrm{Hom}
\left( H^{0}\left( E\right) ,H^{1}\left( E\right) \right) }{image\left(
H^{1}\left( EndE\right) \rightarrow \mathrm{Hom}\left( H^{0}\left( E\right)
,H^{1}\left( E\right) \right) \right) }
\end{equation*}
on.
\end{proof}

\section{The BFM-conjecture for globally generated bundles}

In this final section we prove the injectivity of the map 
\begin{equation*}
\mu _{0}:H^{0}\left( E\right) \otimes H^{0}\left( E\right) \rightarrow
H^{0}\left( E\otimes E\right)
\end{equation*}
in the case of a generic curve $C$ of genus $g>1$ and a vector bundle $E$
satisfying Condition 1 and Condition 2 of \S 1. (Notice that we have only
established the existence of such $E$ when $g>3$ since we relied on the
existence of a non-hyperelliptic $E$ with one vanishing theta-null.)

The isomorphism $\left( \ref{Niso}\right) $ and Serre duality gives 
\begin{equation*}
H^{1}\left( E\right) =H^{1}\left( E^{\vee }\otimes \omega \right)
=H^{0}\left( E\right) ^{\vee }
\end{equation*}
so that we have a decompostion of 
\begin{equation*}
\mathrm{Hom}\left( H^{0}\left( E\right) ,H^{1}\left( E\right) \right)
=H^{1}\left( E\right) \otimes H^{1}\left( E\right)
\end{equation*}
into the direct sum 
\begin{equation*}
Sym^{2}H^{1}\left( E\right) \oplus \bigwedge\nolimits^{2}H^{1}\left(
E\right) .
\end{equation*}
Looking at the $\pm 1$-eigenspaces of the natural action given by reversing
factors, the map 
\begin{equation*}
H^{1}\left( EndE\right) \rightarrow \mathrm{Hom}\left( H^{0}\left( E\right)
,H^{1}\left( E\right) \right)
\end{equation*}
is the direct sum of the maps 
\begin{equation*}
H^{1}\left( End^{0}E\right) \rightarrow Sym^{2}H^{1}\left( E\right)
\end{equation*}
and 
\begin{equation*}
H^{1}\left( \mathcal{O}_{C}\right) \rightarrow
\bigwedge\nolimits^{2}H^{1}\left( E\right) .
\end{equation*}
Dually the multiplication map $\mu _{0}$ is the direct sum of the maps 
\begin{equation*}
Sym^{2}\left( H^{0}\left( E\right) \right) \rightarrow H^{0}\left(
Sym^{2}E\right)
\end{equation*}
and 
\begin{equation*}
\bigwedge\nolimits^{2}H^{0}\left( E\right) \rightarrow H^{0}\left(
\bigwedge\nolimits^{2}E\right) =H^{0}\left( \omega \right) .
\end{equation*}

This final argument generalizes a similar one in \S 9.14 of \cite{ACGH}. Let 
\begin{equation*}
p_{i}:C\times C\rightarrow C
\end{equation*}
for $i=1,2$ denote the two projections and let $D\subseteq C\times C$ denote
the diagonal and $\mathcal{I}$ its ideal in $\mathcal{O}_{C\times C}$. We
have 
\begin{eqnarray*}
\frak{D}_{n}\left( E\right) &=&Hom\left( \left( p_{1}\right) _{*}\frac{
p_{2}^{*}E}{\mathcal{I}^{n+1}p_{2}^{*}E},E\right) \\
&=&Hom\left( \left( p_{1}\right) _{*}\frac{p_{2}^{*}E}{\mathcal{I}
^{n+1}p_{2}^{*}E},E^{\vee }\otimes \omega \right) \\
&=&Hom\left( \left( p_{1}\right) _{*}\frac{E\boxtimes E}{\mathcal{I}
^{n+1}\left( E\boxtimes E\right) },\omega \right) .
\end{eqnarray*}
This gives a perfect pairing 
\begin{equation*}
\frak{D}_{n}\left( E\right) \otimes \left( p_{1}\right) _{*}\frac{E\boxtimes
E}{\mathcal{I}^{n+1}\left( E\boxtimes E\right) }\rightarrow \omega
\end{equation*}
under which the exact sequence 
\begin{equation*}
0\rightarrow \frak{D}_{n-1}\left( E\right) \rightarrow \frak{D}_{n}\left(
E\right) \overset{\sigma }{\longrightarrow }T_{C}^{n}\otimes EndE\rightarrow
0
\end{equation*}
pairs with the exact sequence 
\begin{equation*}
0\rightarrow \left( p_{1}\right) _{*}\frac{\mathcal{I}^{n}E\boxtimes E}{%
\mathcal{I}^{n+1}\left( E\boxtimes E\right) }\overset{\sigma ^{\vee }}{
\longrightarrow }\left( p_{1}\right) _{*}\frac{E\boxtimes E}{\mathcal{I}%
^{n+1}\left( E\boxtimes E\right) }\rightarrow \left( p_{1}\right) _{*}\frac{
E\boxtimes E}{\mathcal{I}^{n}\left( E\boxtimes E\right) }\rightarrow 0
\end{equation*}
where 
\begin{equation*}
\left( p_{1}\right) _{*}\frac{\mathcal{I}^{n}\left( E\boxtimes E\right) }{%
\mathcal{I}^{n+1}\left( E\boxtimes E\right) }=\left( \left( p_{1}\right) _{*}%
\frac{\mathcal{I}^{n}}{\mathcal{I}^{n+1}}\right) \otimes E\otimes E.
\end{equation*}
So we have induced a perfect pairing 
\begin{equation*}
\frac{\frak{D}_{n}\left( E\right) }{\frak{D}_{n-1}\left( E\right) }\otimes
\left( \left( p_{1}\right) _{*}\frac{\mathcal{I}^{n}}{\mathcal{I}^{n+1}}
\otimes E\otimes E\right) \rightarrow \omega .
\end{equation*}

Now the natural maps 
\begin{equation*}
\tilde{\mu}_{n}:H^{0}\left( E\right) \otimes H^{0}\left( E\right)
=H^{0}\left( E\boxtimes E\right) \rightarrow H^{0}\left( \left( p_{1}\right)
_{*}\frac{E\boxtimes E}{\mathcal{I}^{n+1}\left( E\boxtimes E\right) }\right)
\end{equation*}
are the adjoints of the maps 
\begin{equation*}
\tilde{\mu}^{n}:H^{1}\left( \frak{D}_{n}\left( E\right) \right) \rightarrow 
\mathrm{Hom}\left( H^{0}\left( E\right) ,H^{1}\left( E\right) \right)
=H^{1}\left( E\right) \otimes H^{1}\left( E\right)
\end{equation*}
and the induced maps 
\begin{equation*}
\mu _{n}:\ker \tilde{\mu}_{n-1}=H^{0}\left( \left( p_{1}\right) _{*}\mathcal{%
\ I}^{n}\left( E\boxtimes E\right) \right) \rightarrow H^{0}\left( \left(
p_{1}\right) _{*}\frac{\mathcal{I}^{n}\left( E\boxtimes E\right) }{\mathcal{I%
}^{n+1}\left( E\boxtimes E\right) }\right)
\end{equation*}
are the adjoints of the maps 
\begin{equation*}
\mu ^{n}:H^{1}\left( T_{C}^{n}\otimes EndE\right) \rightarrow \frac{\mathrm{%
Hom}\left( H^{0}\left( E\right) ,H^{1}\left( E\right) \right) }{image\tilde{%
\mu}^{n-1}}.
\end{equation*}
Thus the vanishing of $\tilde{\mu}^{n}$ for $n\geq 1$ implies the vanishing
of the $\mu _{n}$ for $n\geq 1$. But this in turn implies that the maps 
\begin{equation*}
H^{0}\left( \mathcal{I}^{n+1}\left( E\boxtimes E\right) \right) \rightarrow
H^{0}\left( \mathcal{I}^{n}\left( E\boxtimes E\right) \right)
\end{equation*}
are isomorphisms for all $n\geq 1$. Since 
\begin{equation*}
h^{0}\left( \mathcal{I}^{n}\left( E\boxtimes E\right) \right) =0
\end{equation*}
for large $n$ we conclude that 
\begin{equation*}
\ker \tilde{\mu}_{0}=0.
\end{equation*}
But 
\begin{equation*}
\tilde{\mu}_{0}=\mu _{0}.
\end{equation*}

\section{Appendix: Mercat's example}

\subsection{Necessity of Condition 2}

Theorem \ref{MT}, even injectivity on the symmetric summand, is definitely
false without Condition 2. We are indebted to Vincent Mercat for pointing
out a counterexample. Suppose $g$ is odd and $\Bbb{P}\left( H^{0}\left(
B\right) \right) $ is a general (basepoint-free) $g_{\frac{g+3}{2}}^{1}$'s
on a general curve $C$. Since $C$ is generic, the injectivity of the Petri
map 
\begin{equation*}
H^{0}\left( B^{2}\right) \otimes H^{0}\left( \omega \otimes B^{-2}\right)
\rightarrow H^{0}\left( \omega \right)
\end{equation*}
and the injectivity of 
\begin{equation*}
Sym^{2}H^{0}\left( B\right) \rightarrow H^{0}\left( B^{2}\right)
\end{equation*}
imply that $h^{0}\left( \omega \otimes B^{-2}\right) =0$. So 
\begin{equation}
\dim Ext^{1}\left( B,\omega \otimes B^{-1}\right) =4  \label{ext}
\end{equation}
so that, up to non-zero scalar, there is a unique non-trivial $\varepsilon
\in Ext^{1}\left( B,\omega \otimes B^{-1}\right) $ such that the (symmetric)
map 
\begin{equation*}
\varepsilon \cdot :H^{0}\left( B\right) \rightarrow H^{1}\left( \omega
\otimes B^{-1}\right)
\end{equation*}
is zero. In fact $\varepsilon $ is a hyperplane in $H^{0}\left( B^{2}\right) 
$ given by the image of the injective map 
\begin{equation*}
Sym^{2}H^{0}\left( B\right) \rightarrow H^{0}\left( B^{2}\right) ,
\end{equation*}
then $\varepsilon $ gives an extension 
\begin{equation}
0\rightarrow \omega \otimes B^{-1}\rightarrow E\rightarrow B\rightarrow 0
\label{anotherext}
\end{equation}
such that the induced map 
\begin{equation*}
H^{0}\left( E\right) \rightarrow H^{0}\left( B\right)
\end{equation*}
is surjective. Thus $E$ is globally generated. By degree, every endomorphism
of $E$ acts as multiplication by a scalar on $\omega \otimes B^{-1}$ and it
follows easily that 
\begin{equation*}
h^{0}\left( End^{0}E\right) =0,
\end{equation*}
that is, $E$ is simple. Also, since $h^{0}\left( \omega \otimes
B^{-1}\right) =\frac{g-1}{2}$, 
\begin{equation*}
h^{0}\left( E\right) =\frac{g+3}{2}.
\end{equation*}
So, for large enough $g$, 
\begin{equation*}
\dim Sym^{2}H^{0}\left( E\right) >h^{0}\left( Sym^{2}E\right) .
\end{equation*}

\subsection{Geometry of the spectral curve in Mercat's example}

In the presentation 
\begin{equation*}
0\rightarrow F\rightarrow H^{0}\left( E\right) \otimes \mathcal{O}_{C}%
\overset{\varepsilon }{\longrightarrow }E\rightarrow 0
\end{equation*}
in Mercat's example, the sections of $E$ coming from sections of $\omega
\otimes B^{-1}$ all have 
\begin{equation*}
2g-2-\frac{g+3}{2}=\frac{3g-7}{2}
\end{equation*}
zeros on $C$. For 
\begin{equation*}
W=H^{0}\left( E\right)
\end{equation*}
the subspace 
\begin{equation*}
W^{\prime }:=H^{0}\left( \omega \otimes B^{-1}\right) \subseteq H^{0}\left(
E\right)
\end{equation*}
gives rise to an exact sequence 
\begin{equation*}
0\rightarrow F^{\prime }\rightarrow W^{\prime }\otimes \mathcal{O}%
_{C}\rightarrow \omega \otimes B^{-1}\rightarrow 0
\end{equation*}
and so a commutative diagram 
\begin{equation}
\begin{array}{ccccccccc}
&  & 0 &  & 0 &  & 0 &  &  \\ 
&  & \downarrow &  & \downarrow &  & \downarrow &  &  \\ 
0 & \rightarrow & B^{\vee } & \rightarrow & \left( \frac{W}{W^{\prime }}%
\right) ^{\vee }\otimes \mathcal{O}_{C} & \rightarrow & B & \rightarrow & 0
\\ 
&  & \downarrow &  & \downarrow &  & \downarrow &  &  \\ 
0 & \rightarrow & E^{\vee } & \rightarrow & W^{\vee }\otimes \mathcal{O}_{C}
& \rightarrow & F^{\vee } & \rightarrow & 0 \\ 
&  & \downarrow &  & \downarrow &  & \downarrow &  &  \\ 
0 & \rightarrow & \omega ^{-1}\otimes B & \rightarrow & \left( W^{\prime
}\right) ^{\vee }\otimes \mathcal{O}_{C} & \rightarrow & \left( F^{\prime
}\right) ^{\vee } & \rightarrow & 0 \\ 
&  & \downarrow &  & \downarrow &  & \downarrow &  &  \\ 
&  & 0 &  & 0 &  & 0 &  & 
\end{array}
\label{diag}
\end{equation}
with exact rows and columns. The sheaf of linear functionals on $E$ which
vanish on $\omega \otimes B^{-1}$ lie in the intersection of $\left( \frac{W%
}{W^{\prime }}\right) ^{\vee }\otimes \mathcal{O}_{C}$ and $E^{\vee }$ in $%
W^{\vee }\otimes \mathcal{O}_{C}$. Thus, for the mapping 
\begin{eqnarray*}
p &:&C\rightarrow G, \\
c &\mapsto &\Bbb{P}\left( \varepsilon ^{\vee }\left( E_{c}^{\vee }\right)
\right)
\end{eqnarray*}
$p\left( C\right) $ lies in the Schubert cycle of lines incident to the line 
\begin{equation*}
l_{0}:=\Bbb{P}\left( \left( \frac{W}{W^{\prime }}\right) ^{\vee }\right)
\subseteq \Bbb{P}\left( W^{\vee }\right) .
\end{equation*}
In fact, from $\left( \ref{diag}\right) $ it is easy to see that the $g_{%
\frac{g+3}{2}}^{1}$ on $C$ is given by assigning to each point $x\in l_{0}$
the divisor of points $c\in C$ such that $x\in \Bbb{P}\left( E_{c}^{\vee
}\right) $.

Projection with center $l_{0}$ maps 
\begin{equation*}
\tilde{C}\subseteq \Bbb{P}\left( W^{\vee }\right)
\end{equation*}
$2-1$ to the embedding 
\begin{equation*}
C\subseteq \Bbb{P}\left( \left( W^{\prime }\right) ^{\vee }\right)
\end{equation*}
induced by the bundle inclusion 
\begin{equation*}
\omega ^{-1}\otimes B\subseteq \left( W^{\prime }\right) ^{\vee }\otimes 
\mathcal{O}_{C},
\end{equation*}
that is, by the sections of the line bundle $\omega \otimes B^{-1}$. Also,
via intersection, the line $l_{0}$ gives a section 
\begin{equation*}
S:=\Bbb{P}\left( B^{\vee }\right)
\end{equation*}
of the bundle 
\begin{equation*}
\pi :\Bbb{P}\left( E^{\vee }\right) \rightarrow C
\end{equation*}
so that $\left( \ref{anotherext}\right) $ is obtained by applying 
\begin{equation*}
\pi _{*}\circ \left( \mathcal{O}_{\Bbb{P}\left( E^{\vee }\right) }\left(
1\right) \otimes \ \right)
\end{equation*}
to the sequence 
\begin{equation*}
0\rightarrow \mathcal{O}_{\Bbb{P}\left( E^{\vee }\right) }\left( -S\right)
\rightarrow \mathcal{O}_{\Bbb{P}\left( E^{\vee }\right) }\rightarrow 
\mathcal{O}_{S}\rightarrow 0.
\end{equation*}

We next choose a general 
\begin{equation*}
\tilde{\varphi}\in Hom^{0}\left( E^{\vee },E\right)
\end{equation*}
with spectral double cover 
\begin{equation*}
\tilde{C}\subseteq \Bbb{P}\left( E^{\vee }\right) \rightarrow C.
\end{equation*}
For the line bundle 
\begin{equation*}
L=\left. \mathcal{O}_{\Bbb{P}\left( E^{\vee }\right) }\left( 1\right)
\right| _{\tilde{C}}
\end{equation*}
and involution 
\begin{equation*}
\iota :\tilde{C}\rightarrow \tilde{C},
\end{equation*}
we have 
\begin{equation}
\pi _{*}\left( \iota ^{*}L\right) =\pi _{*}L.  \label{Eiso}
\end{equation}

The inclusion 
\begin{equation*}
\omega \otimes B^{-1}\rightarrow \pi _{*}L
\end{equation*}
induces a non-trivial map 
\begin{equation}
\pi ^{*}\left( \omega \otimes B^{-1}\right) \rightarrow L  \label{incl}
\end{equation}
so that 
\begin{eqnarray*}
L &=&\pi ^{*}\left( \omega \otimes B^{-1}\right) \otimes \mathcal{O}\left( 
\tilde{D}_{0}\right) \\
\iota ^{*}L &=&\pi ^{*}\left( \omega \otimes B^{-1}\right) \otimes \mathcal{O%
}\left( \iota ^{*}\tilde{D}_{0}\right)
\end{eqnarray*}
for some effective divisor $\tilde{D}_{0}$ on $\tilde{C}$ with 
\begin{eqnarray*}
\omega &=&\det \pi _{*}L=\omega \otimes B^{-2}\otimes \mathcal{O}\left( \pi
_{*}\tilde{D}_{0}\right) \\
\mathcal{O}\left( \pi _{*}\tilde{D}_{0}\right) &=&B^{2}.
\end{eqnarray*}
In fact, since 
\begin{eqnarray*}
\left. \mathcal{O}_{\Bbb{P}\left( W^{\vee }\right) }\left( 1\right) \right|
_{\tilde{C}} &=&\pi ^{*}\left( \omega \otimes B^{-1}\right) \otimes \mathcal{%
\ O}\left( \tilde{D}_{0}\right) \\
\left. \mathcal{O}_{\Bbb{P}\left( \left( W^{\prime }\right) ^{\vee }\right)
}\left( 1\right) \right| _{C} &=&\omega \otimes B^{-1}
\end{eqnarray*}
the support of $\tilde{D}_{0}$ must equal the intersection of $\varepsilon
^{\vee }\left( \tilde{C}\right) $ with the center $l_{0}$ of the projection 
\begin{equation*}
\Bbb{P}\left( W^{\vee }\right) \dashrightarrow \Bbb{P}\left( \left(
W^{\prime }\right) ^{\vee }\right) .
\end{equation*}

Now suppose, in Mercat's example, 
\begin{equation*}
\tilde{\varphi}=\tilde{Q}
\end{equation*}
for some quadric 
\begin{equation*}
Q\in Sym^{2}\left( H^{0}\left( E\right) \right) .
\end{equation*}
For example, suppose that 
\begin{equation*}
Sym^{2}H^{0}\left( \mathcal{O}_{\Bbb{P}\left( E^{\vee }\right) }\left(
1\right) \right) =Sym^{2}H^{0}\left( E\right) \rightarrow H^{0}\left(
Sym^{2}E\right) =H^{0}\left( \mathcal{O}_{\Bbb{P}\left( E^{\vee }\right)
}\left( 2\right) \right)
\end{equation*}
is surjective. Since 
\begin{equation*}
\varepsilon ^{\vee }\left( \tilde{C}\right) =Q\cap \varepsilon ^{\vee
}\left( E^{\vee }\right) ,
\end{equation*}
and the $g_{\frac{g+3}{2}}^{1}$ on $C$ is given by assigning to each point $%
x\in l_{0}$ the divisor of points $c\in C$ such that $x\in \Bbb{P}\left(
E_{c}^{\vee }\right) $, we conclude that 
\begin{equation*}
\tilde{D}_{0}=\tilde{D}_{1}+\tilde{D}_{2},
\end{equation*}
where the divisors 
\begin{equation*}
\pi _{*}\tilde{D}_{1},\pi _{*}\tilde{D}_{2}\in g_{\frac{g+3}{2}}^{1}
\end{equation*}
are those parametrized by the two points of 
\begin{equation*}
l_{0}\cap Q.
\end{equation*}
Then 
\begin{equation*}
H^{0}\left( \mathcal{O}_{\Bbb{P}\left( E^{\vee }\right) }\left( 2\right)
\right) \rightarrow H^{0}\left( \mathcal{O}_{S}\left( 2\right) \right)
=H^{0}\left( B^{2}\right)
\end{equation*}
is not surjective since it has image $Sym^{2}H^{0}\left( B\right)
\subsetneqq H^{0}\left( B^{2}\right) $.

By \S 1.3.2 of \cite[B]{B} Condition 2 is equivalent to the condition 
\begin{equation*}
h^{0}\left( \pi ^{*}\omega \otimes \mathcal{O}\left( \iota ^{*}\tilde{D}_{0}-%
\tilde{D}_{0}\right) \right) =0.
\end{equation*}
It follows from the main theorem of this paper that Condition 2 must fail
for Mercat's example. It would be reassuring to have an independent proof of
this fact, but the authors have not found an independent method for deciding
whether or not Condition 2 holds in the case of Mercat's example.

\end{document}